\documentstyle[11pt,a4wide]{article}
\font \tenmsb=msbm10 scaled \magstep 1
\font \sevenmsb=msbm7 scaled \magstep 1
\font \fivemsb=msbm5 scaled \magstep 1
\newfam \msbfam
\textfont \msbfam=\tenmsb
\scriptfont \msbfam=\sevenmsb
\scriptscriptfont \msbfam=\fivemsb
\def \Bbb#1{\fam \msbfam \relax#1}

\font\teneufm=eufm10 scaled \magstep 1
\font\seveneufm=eufm7 scaled \magstep 1
\font\fiveeufm=eufm5 scaled \magstep 1
\newfam\eufmfam
\textfont\eufmfam=\teneufm
\scriptfont\eufmfam=\seveneufm
\scriptscriptfont\eufmfam=\fiveeufm
\def\frak#1{{\fam\eufmfam\relax#1}}

\title{\bf FUNCTION THEORY IN THE QUANTUM MATRIX BALL:
 AN INVARIANT INTEGRAL}

\author{\sl D. Shklyarov$^\ddagger$ \and \sl S. Sinel'shchikov$^\dagger$
\and \sl L. Vaksman$^\dagger$}

\date{\tt $^\dagger$Institute for Low Temperature Physics \& Engineering\\
National Academy of Sciences of Ukraine\\ \ \\
$^\ddagger$Kharkov State University\\
Department of Mechanics and Mathematics}

\begin{document}

\maketitle

 {\bf 1.} It is well known \cite{H} that the unit ball $U$  in the space of
rectangle complex matrices is a bounded symmetric domain. Our recent work
\cite{SV} introduces q-analogues of all such domains. This work considers
the q-analogue of the matrix ball and presents an explicit formula for a
positive invariant integral.

 A parameter $q$ involved into the formulations of the main results is
assumed to be a number: $0<q<1$. However, in sections 2, 3 the ground field
will be not ${\Bbb C}$ but ${\Bbb C}(q^{1/2})$, the field of rational
functions of the parameter $q^{1/2}$ (cf. \cite{J}).

\bigskip

 {\bf 2.} The Hopf algebra $U_q \frak{sl}_N$ is determined by its generators
$\{E_i,F_i,K_i^{\pm 1}\}_{i=1,\ldots,N-1}$ and the well known Drinfeld-Jimbo
relations (see \cite{J}). The universal enveloping algebra $U \frak{sl}_N$
could be derived from $U_q \frak{sl}_N$ via the "change of variables"
$q=e^{-h/2}$, $K_j^{\pm 1}=e^{\mp hH_j/2}$, $j=1,\ldots,N-1$, and a formal
passage to a limit as $h \to 0$ in the list of relations.

 Everywhere in the sequel $m,n \in{\Bbb N}$, $N=m+n$. We use the standard
notation $\frak{su}_{nm}$ for the Lie algebra of the automorphism group of
the unit matrix ball $U \subset{\rm Mat}_{mn}$.

 Equip the Hopf algebra $U_q \frak{sl}_N$ with the involution, which is
defined on the generators $K_j^{\pm 1}$, $E_j$, $F_j$, $j=1,\ldots,N-1$ by
$$\left(K_j^{\pm 1}\right)^*=K_j^{\pm 1},\qquad E_j^*=\left \{
\begin{array}{rr}K_jF_j,& j \ne n \\ -K_jF_j,&j=n \end{array}\right.,
\qquad F_j^*=\left \{\begin{array}{rr}E_jK_j^{-1},& j \ne n \\
-E_jK_j^{-1},&j=n \end{array}\right..$$
The Hopf $*$-algebra $U_q \frak{su}_{nm}=(U_q \frak{sl}_N,*)$ arising in
    this way is a q-analogue of the Hopf algebra $U \frak{su}_{nm}$.

 Remind some well known definitions $\cite{CP}$. An algebra $F$ is said to
be an $A$-module algebra if it is a module over a Hopf algebra $A$  and the
multiplication $F \otimes F \to F$, $f_1 \otimes f_2 \mapsto f_1 \cdot f_2$,
is a morphism of $A$-modules. In the case of a $*$-algebra $F$ and a
Hopf-$*$-algebra $A$, there is an additional requirement that the
involutions agree as follows:
$$(af)^*=(S(a))^*f^*,\qquad a \in A,\;f \in F,$$
with $S:A \to A$ being the antipode of $A$.

\bigskip

 {\bf 3.} In \cite{SV} a $U_q \frak{su}_{nm}$-module algebra ${\rm
Pol}({\rm Mat}_{mn})_q$ and its $U_q \frak{sl}_N$-module subalgebra ${\Bbb
C}[{\rm Mat}_{mn}]_q$ were introduced (the notation ${\frak g}_{-1}$ was
used in \cite{SV} instead of ${\rm Mat}_{mn}$). These algebras are
q-analogues of polynomial algebras in the vector spaces ${\rm
Mat}_{mn}$. For that, the passage to dual coalgebras was implemented,
which constitutes an approach of V.Drinfeld \cite{D}. We present below a
description of the algebra ${\Bbb C}[{\rm Mat}_{mn}]_q$ in terms of
generators and relations.

 With the definitions of \cite{SV} as a background, one can prove the
following two propositions.

\medskip

 {\sc Proposition 1.} {\it There exists a unique family $\{z_a^\alpha
\}_{a=1,\ldots,n;\alpha=1,\ldots,m}$ of elements of the $U_q
\frak{sl}_N$-module algebra ${\Bbb C}[{\rm Mat}_{mn}]_q$ such that for all
${a=1,\ldots,n;\alpha=1,\ldots,m}$
$$H_nz_a^\alpha=\left \{\begin{array}{ccl}2z_a^\alpha &,&a=n \;\&\;\alpha=m
\\ z_a^\alpha &,&a=n \;\&\;\alpha \ne m \quad{\rm or}\quad a \ne n \;\&\;
\alpha=m \\ 0 &,&\;{\rm otherwise}\end{array}\right.\eqno(1)$$
$$F_nz_a^\alpha=q^{1/2}\cdot \left \{\begin{array}{ccl}1 &,& a=n \;\&
\;\alpha=m \\ 0 &,&\;{\rm otherwise}\end{array}\right.\eqno(2)$$
$$E_nz_a^\alpha=-q^{1/2}\cdot \left \{\begin{array}{ccl}q^{-1}z_a^mz_n^\alpha
 &,&a \ne n \;\&\;\alpha \ne m \\
(z_n^m)^2 &,& a=n \;\&\;\alpha=m \\ z_n^mz_a^{\alpha} &,&\;{\rm
otherwise}\end{array}\right.\eqno(3)$$
and with $k \ne n$
$$H_kz_a^\alpha= \left \{\begin{array}{ccl}z_a^\alpha &,& k<n
\;\&\;a=k \quad{\rm or}\quad k>n \;\&\;\alpha=N-k \\-z_a^\alpha &,& k<n
\;\&\;a=k+1 \quad{\rm or}\quad k>n \;\&\;\alpha=N-k+1 \\ 0 &,&\;{\rm
otherwise}\end{array}\right.$$
$$F_kz_a^\alpha=q^{1/2}\cdot \left \{\begin{array}{ccl}z_{a+1}^\alpha &,&
k<n \;\&\;a=k \\ z_a^{\alpha+1} &,& k>n \;\&\;\alpha=N-k \\ 0 &,&\;{\rm
otherwise}\end{array}\right.$$
$$E_kz_a^\alpha=q^{-1/2}\cdot \left \{\begin{array}{ccl}z_{a-1}^\alpha &,&
k<n \;\&\; a=k+1 \\ z_a^{\alpha-1}&,& k>n \;\&\;\alpha=N-k+1 \\ 0 &,&
\;{\rm otherwise}\end{array}\right.$$}

\medskip

 {\sc Proposition 2.} {\it $\{z_a^\alpha
\}_{a=1,\ldots,n;\alpha=1,\ldots,m}$ generate ${\Bbb C}[{\rm Mat}_{mn}]_q$
as an algebra and ${\rm Pol}({\rm Mat}_{mn})_q$ as a $*$-algebra. The
complete list of relations is  following:
$$z_a^\alpha z_b^\beta=\left \{\begin{array}{ccl}qz_b^\beta z_a^\alpha &,&
a=b \;\&\;\alpha<\beta \quad{\rm or}\quad a<b \;\&\;\alpha=\beta \\
z_b^\beta z_a^\alpha &,& a<b \;\&\;\alpha>\beta \\ z_b^\beta
z_a^\alpha+(q-q^{-1})z_a^\beta z_b^\alpha &,& a<b \;\&\;\alpha<\beta
\end{array} \right.$$
$$(z_b^\beta)^*z_a^\alpha=q^2 \cdot \sum_{a^\prime,b^\prime=1}^n
\sum_{\alpha^\prime,\beta^\prime=1}^m R_{ba}^{\prime b^\prime
a^\prime}R_{\beta^\prime \alpha^\prime}^{\prime \prime \beta
\alpha}\cdot
z_{a^\prime}^{\alpha^\prime}\left(z_{b^\prime}^{\beta^\prime}\right)^*+
(1-q^2)\delta_{ab}\delta^{\alpha \beta},$$
with $\delta_{ab}$, $\delta^{\alpha \beta}$ being the Kronecker symbols and
$$R_{ba}^{\prime b^\prime a^\prime}=\left \{\begin{array}{ccl}q^{-1} &,& a
\ne b \;\&\;b=b^\prime \;\&\;a=a^\prime \\ 1 &,& a=b=a^\prime=b^\prime \\
-(q^{-2}-1) &,& a=b \;\&\;a^\prime=b^\prime \;\&\;a^\prime>a \\ 0 &,&\;{\rm
otherwise}\end{array}\right.$$
$$R_{\beta^\prime \alpha^\prime}^{\prime \prime \beta \alpha}=\left
\{\begin{array}{ccl}q^{-1} &,& \alpha \ne \beta \;\&\; \beta=\beta^\prime
\;\&\;\alpha=\alpha^\prime \\ 1 &,& \alpha=\beta=\alpha^\prime=\beta^\prime
\\-(q^{-2}-1) &,& \alpha=\beta \;\&\; \alpha^\prime=\beta^\prime
\;\&\;\alpha^\prime>\alpha \\ 0 &,& \;{\rm otherwise}\end{array}\right.$$}

\medskip

{\sc Corollary 3.}
$$\left(z_n^m \right)^*z_n^m=q^2z_n^m \left(z_n^m \right)^*+1-q^2.\eqno(4)$$

\medskip

 Note that $z_n^m$ generates a $U_q \frak{su}_{11}$-module algebra ${\rm
Pol}({\Bbb C})_q$ determined by relations (1) -- (4) (see \cite{SV}).
 Commutation relations similar to those given in proposition 2 appear
in a different context in \cite{CHZ}.
\bigskip

 {\bf 4.} Everywhere in the sequel we assume $q \in(0,1)$, and ${\Bbb C}$ is
considered as a ground field. We keep the notation $U_q \frak{su}_{nm}$,
${\rm Pol}({\rm Mat}_{mn})_q$ for the Hopf $*$-algebra and the covariant
$*$-algebra determined by the generators $\{E_j,F_j,K_j^{\pm
1}\}_{j=1,\ldots,N-1}$ and $\{z_a^\alpha
\}_{a=1,\ldots,n;\alpha=1,\ldots,m}$ respectively and the relations as
above.

 It is well known that in the classical case $q=1$ the positive
$SU_{nm}$-invariant measure on the matrix ball is infinite. Thus the
positive invariant integral could not be defined on the polynomial
algebra.In the quantum case this obstacle is still in effect, therefore we
need to extend the $U_q \frak{su}_{nm}$-module algebra ${\rm Pol}({\rm
Mat}_{mn})_q$ up to the $U_q \frak{su}_{nm}$-module algebra ${\rm
Fun}(U)_q$. The construction we produce below can seem non-substantiated. We
refer the reader to the work \cite{SSV2}, where it is described in more
details in an important special case of the quantum disc (i.e. m=n=1).

 Consider the $*$-algebra ${\rm Fun}(U)_q \supset{\rm Pol}({\rm Mat}_{mn})_q$
derived from ${\rm Pol}({\rm Mat}_{mn})_q$ by adjunction an additional
generator $f_0$ such that
$$f_0=f_0^2=f_0^*\quad{\rm and}\quad \left(z_a^\alpha
\right)^*f_0=f_0z_a^\alpha=0, \quad
a=1,\ldots,n;\;\alpha=1,\ldots,m.\eqno(5)$$
((5) allows one to treat $f_0$ as a q-analogue of the delta-function at
zero.)

 Extend the structure of a $U_q \frak{su}_{nm}$-module algebra from ${\rm
Pol}({\rm Mat}_{mn})_q$ onto ${\rm Fun}(U)_q$ via the relations (6)
(which were obtained in our earlier work \cite{SSV2} in the special case
$m=n=1$) and (7). The following proposition is a consequence of the 
definition of the
$*$-algebra ${\rm Fun}(U)_q$, the relation (4), and the quasicommutativity
of $z_n^m$ with $\left(z_a^\alpha \right)^*$ for $(n,m)\ne(a,\alpha)$.

\medskip

 {\sc Proposition 4.} {\it There exists a unique extension of the structure
of a $U_q \frak{su}_{nm}$-module algebra from ${\rm Pol}({\rm Mat}_{mn})_q$
onto ${\rm Fun}(U)_q$ such that
$$H_nf_0=0,\qquad F_nf_0=-{q^{1/2}\over q^{-2}-1}f_0 \cdot \left(z_n^m
\right)^*,\qquad E_nf_0=-{q^{1/2}\over 1-q^2}z_n^m \cdot f_0 \eqno(6)$$
and with $k \ne n$
$$H_kf_0=F_kf_0=E_kf_0=0.\eqno(7)$$}

\medskip

 {\bf 5.} The two-sided ideal $D(U)_q \stackrel{\rm def}{=}{\rm
Fun}(U)_qf_0{\rm Fun}(U)_q$ is a $U_q \frak{su}_{nm}$-module algebra. Its
elements will be called the finite functions in the quantum matrix ball.

  Our purpose is to produce a positive invariant integral on the
algebra of finite functions
$$D(U)_q \to{\Bbb C},\qquad f \mapsto \int \limits_{U_q}fd \nu,$$
i.e.  such linear functional that $\displaystyle \int 
\limits_{U_q}f^*fd\nu>0$
for all $f \ne 0$, and $\displaystyle \int \limits_{U_q}(\xi f)d
\nu=\varepsilon(\xi)\int \limits_{U_q}fd \nu$ for all $\xi \in U_q
\frak{sl}_N$, $f \in D(U)_q$. Here $\varepsilon$ is the counit of the Hopf
algebra $U_q \frak{sl}_N$.

 Consider the vector space ${\cal H}={\Bbb C}[{\rm Mat}_{mn}]_qf_0 \subset
D(U)_q$. Equip it with the grading (see \cite{SV}):
$${\cal H}=\oplus_{j=0}^\infty{\cal H}_j,\qquad{\cal H}_j \stackrel{\rm
def}{=}\{v \in{\cal H}|\;H_0v=2jv \}$$
with
$$H_0={2 \over m+n}\left(m \sum_{j=1}^{n-1}jH_j+n
\sum_{j=1}^{m-1}jH_{N-j}+mnH_n \right).$$
Evidently, ${\rm dim}\;{\cal H}_j<\infty$ for all $j \in{\Bbb Z}_+$.

 Let $T$ be the representation of $D(U)_q$ in ${\cal H}$ given by $T_f
\psi=f \cdot \psi$, $f \in D(U)_q$, $\psi \in{\cal H}\subset D(U)_q$.
The following statements are direct consequences of the definitions.

\medskip

 {\sc Lemma 5.} {\it For any element $f \in D(U)_q$ there exists such
positive integer $M(f)$ that $T_f{\cal H}_j=0$ for all $j \ge M(f)$}.

\medskip

 {\sc Corollary 6.} {\it All the operators $T_f$, $f \in D(U)_q$, are finite
dimensional.}

\medskip

 It is easy to prove the existence and uniqueness of such scalar product in
${\cal H}$ that $(f_0,f_0)=1$ and $(T_f
\psi_1,\psi_2)=(\psi_1,T_{f^*}\psi_2)$ for all $f \in D(U)_q$,
$\psi_1,\psi_2 \in{\cal H}$.

 The representation $T$ is faithful ($T(f)\ne 0$ for $f \ne 0$), and one
also has

\medskip

  {\it $(\psi,\psi)>0$ for all $\psi \in{\cal H}$, $\psi \ne
0$.}

\bigskip

 {\bf 6.} Let $U_q \frak{p}_+\subset U_q \frak{sl}_N$ be the Hopf subalgebra
generated by all the generators $E_j,F_j,K_j^{\pm 1}$, $j=1,\ldots,N-1$
except $F_n$.

 (6), (7) imply that the subspace ${\cal H}$ is a submodule of the $U_q
\frak{p}_+$-module $D(U)_q$. Let $\Gamma$ stand for the representation of
$U_q \frak{p}_+$ in ${\cal H}$. Let also $\check{\rho}={1 \over 2}\sum
\limits_{j=1}^{N-1}j(N-j)H_j$. Then $e^{h \check{\rho}}=\prod
\limits_{j=1}^{N-1}K_j^{-j(N-j)}$, and the operator $\Gamma(e^{h
\check{\rho}})$ in ${\cal H}$ is well defined and takes each ${\cal H}_j$,
$j \in{\Bbb Z}_+$ into itself. Hence by lemma 5 the trace ${\rm
Tr}(T(f)\Gamma(e^{h \check{\rho}}))$ is well defined for all $f \in D(U)_q$.
An application of theorem 7 yields

\medskip

 {\bf Theorem 8.} {\it The linear functional
$$\int \limits_{U_q}fd \nu \stackrel{\rm def}{=}{\rm Tr}(T(f)\Gamma(e^{h
\check{\rho}}))\eqno(8)$$
on the algebra $D(U)_q$ of finite functions in the quantum matrix ball is
positive and invariant.}

\bigskip

 {\bf 7.} Consider a class of those $U_q \frak{su}_{nm}$-module algebras
for which an analogue of theorem 8 is valid due to obvious reasons. Let
${\cal H}$ be a unitarizable Harish-Chandra module with the lowest weight
over the Hopf algebra $U_q \frak{su}_{nm}$, and let $\Gamma$ be the
associated representation of $U_q \frak{su}_{nm}$ in ${\cal H}$. The
canonical embedding ${\cal H}\otimes{\cal H}^*\hookrightarrow{\rm
End}_{\Bbb C}({\cal H})$ allows one to equip $F={\cal H}\otimes{\cal H}^*$
with a structure of a $U_q \frak{su}_{nm}$-module algebra. Let $T$ be the
associated representation of $F$ in ${\cal H}$. It is well known and can be
easily proved that the linear functional ${\rm Tr}_qT(f)\stackrel{\rm
def}{=}{\rm Tr}(T(f)\Gamma(e^{h \check{\rho}}))$ on $F$ is a positive
invariant integral. It was shown in \cite{SSV5} for the case of quantum disc
that the $U_q \frak{su}_{nm}$-module algebra $D(U)_q$ is a "limit point of
the set of" $U_q \frak{su}_{nm}$-module algebras of the class described
above. This substantiates invariance and positivity of the integral (8).

\end{document}